\documentclass[11pt]{amsart} 
\usepackage{amssymb}

\newcommand{\rightlimit}[1]{\mathop{\lim\limits_{\longrightarrow}}\limits%
                   _{\raise3pt\hbox{$\scriptstyle #1$}}}

\newcommand{\BDelta}{{\bf Sim}}
\newcommand{\bone}{{\mathbf{1}}}
\newcommand{\BZ}{{\mathbb Z}}
\newcommand{\BR}{{\mathbb R}}
\newcommand{\cA}{{\mathcal A}}
\newcommand{\cI}{{\mathcal I}}
\newcommand{\iso}{\buildrel{\sim}\over{\longrightarrow}}
\newcommand{\mono}{\hookrightarrow}
\newcommand{\Ifrom}{{\mathop{I}\limits_{\leftarrow}}}
\newcommand{\N}{N_{\cyc}}
\newcommand{\Sfrom}{{\mathop{S^1}\limits_{\leftarrow}}}
\newcommand{\xy}{{\mathop{xy}\limits^{\longrightarrow}}}
\newcommand{\yx}{{\mathop{yx}\limits^{\longrightarrow}}}
\newcommand{\xz}{{\mathop{xz}\limits^{\longrightarrow}}}
\newcommand{\yz}{{\mathop{yz}\limits^{\longrightarrow}}}
\newcommand{\ZFunct}{\BZ_+\mbox{\rm -Funct}}
\newcommand{\Ztors}{\BZ_+\mbox{\rm -tors}}

\DeclareMathOperator{\Aut}{{Aut}}
\DeclareMathOperator{\byg}{{big}}
\DeclareMathOperator{\cyc}{{cyc}}
\DeclareMathOperator{\Card}{{Card}}
\DeclareMathOperator{\End}{{End}}
\DeclareMathOperator{\Ext}{{Ext}}
\DeclareMathOperator{\Fib}{{Fib}}
\DeclareMathOperator{\Funct}{{Funct}}
\DeclareMathOperator{\Hom}{{Hom}}
\DeclareMathOperator{\hug}{{huge}}
\DeclareMathOperator{\id}{{id}}
\DeclareMathOperator{\Ind}{{Ind-}}
\DeclareMathOperator{\Mor}{{Mor}}
\DeclareMathOperator{\Ob}{{Ob}}
\DeclareMathOperator{\Pro}{{Pro-}}
\DeclareMathOperator{\Sets}{{Sets}}

\newtheorem{pr}{Proposition}[section]

\numberwithin{equation}{section}

\begin{document}

\thanks{Partially supported by NSF grant DMS-0100108}

\title{On the notion of geometric realization}

\author{Vladimir Drinfeld}

\address{Dept. of Math., Univ. of Chicago, 5734 University 
Ave., Chicago, IL 60637, USA}

\email{drinfeld@math.uchicago.edu}

\dedicatory{To Borya Feigin on his 50th birthday}

\begin{abstract}
We explain why geometric realization commutes with Cartesian
products and why the geometric realization of a simplicial
set (resp. cyclic set) is equipped with an action of the
group of orientation preserving homeomorphisms of
the segment $[0,1]$ (resp. the circle).

2000 Math. Subj. Class. 18G30, 55U10, 19D55

\end{abstract}

\maketitle

{{{\leftskip=2.00in {
{{\noindent \it  The next simplest curve is a
circle. Even so simple a figure as this has given rise to so
many and such profound investigations that they could
constitute a course by themselves.}

\smallskip 

\noindent  D.~Hilbert and S.~Cohn-Vossen, Geometry and the
imagination (Anschauliche Geometrie), \linebreak ch.~1, page
1.}}\par}

}}

\bigskip

In this note there are no theorems, its only goal is to
clarify the notion of geometric realization for simplicial
sets \cite{Mi,M,GZ,GJ} and cyclic sets \cite{DHK,BF,L}. We
reformulate the definitions so that the following facts
become obvious:

(i) geometric realization commutes with finite
projective limits (e.g., with Cartesian products);

(ii) the geometric realization of a simplicial set
(resp. cyclic set) is equipped with an action of the group
of orientation preserving homeomorphisms of
the segment $I:=[0,1]$ (resp. the circle $S^1$).

In the traditional approach \cite{Mi,M,GZ,GJ,DHK,BF,L}
these statements are theorems, and understanding their proofs
requires some efforts. 

\medskip

\noindent {\bf Example.} To a small category $C$ there
corresponds a simplicial set $NC$ (the nerve of $C$) and a
cyclic set $\N C$ (the cyclic nerve). It follows from our
formula (\ref{def}) that a point of the geometric
realization $|NC|$  is a {\it piecewise constant functor\,}
$I\to C$ . Here it is assumed that the  category structure
on $I$ comes from the standard order on $I$, and the
definition of piecewise constant functor is explained in
(\ref{nerve}). In \S\ref{3} we give a similar description of
$|\N C|$ in which $I$ is replaced by $S^1$.

\medskip

\noindent {\bf Conventions.} The set of non-negative
integers is denoted by $\BZ_+$. Unless stated otherwise, an
ordered set $\cI$ is always equipped with the category
structure such that the number of morphisms from
$i\in\cI$ to $j\in\cI$ equals $1$ if $i\le j$ and $0$
otherwise.

\medskip

The author thanks D.~Arinkin, A.~Beilinson, D.~Grayson, M.~Kapranov, and
P.~May for  useful remarks and references. The epigraph is
due to R.~Bezrukavnikov.

{\bf Having written this article I learned that a very
similar approach had been developed by A.~Besser \cite{B}
and D.~Grayson \cite{G}. It is used in \cite{G} to treat not
only products of simplicial sets but also Quillen's edgewise
subdivision.}

\section{Simplicial sets} \label{1}

Put $I:=[0,1]$. This is an ordered set and therefore a
category.

Recall that a simplicial set is a functor
$X: \Delta^\circ\to\Sets$, where $\Delta$ is the category
whose objects are the ordered sets $[n]:=\{ 0,\ldots ,n\}$,
$n\in\BZ_+$, and whose morphisms are nondecreasing maps.
One can extend such $X$ to a functor
$\Delta_{\byg}^\circ\to\Sets$, where $\Delta_{\byg}$ is the
category of non-empty finite linearly ordered sets. Such an
extension is unique up to unique isomorphism.

We define the {\it geometric realization} of a
simplicial set $X$ to be the set
\begin{equation} \label{def}
|X|:=\rightlimit{F} X (\pi_0(I\setminus F) ),
\end{equation}
where $F$ runs through the set of all finite subsets of $I$
and $\pi_0(I\setminus F)$ is equipped with the natural order
(so $\pi_0(I\setminus F)\in\Delta_{\byg}$ is a quotient of
$I\setminus F$). At this point $|X|$ is viewed merely as a
set, the definition of the topology on $|X|$ will be
recalled later. One can also rewrite (\ref{def}) as 
\begin{equation} \label{Ifrom}
|X|:=X(\Ifrom ), 
\end{equation}
where $\Ifrom$ is the pro-object of
$\Delta_{\byg}$ which is the projective limit of the objects
$\pi_0(I\setminus F)\in\Delta_{\byg}$ over all finite
subsets $F\subset I$.

It immediately follows from the definition that geometric 
realization commutes with finite projective limits (e.g., the
map $|X\times Y|\to |X|\times |Y|$ is bijective) and that
the group $\Aut I$ of orientation preserving homeomorphisms
of $I$ acts on $|X|$.

\medskip

\noindent {\bf Example.} The geometric realization $|[n]|$ of
$[n]$ (i.e., of the functor $\Delta^{\circ}\to\Sets$
represented by $[n]$) is the set of piecewise constant
nondecreasing functions $f:I\to [n]$ modulo the following
equivalence relation: $f_1\sim f_2$ if $f_1$ and $f_2$ are
equal outside a finite set. We identify $|[n]|$ with the
standard simplex
\begin{equation} \label{sim}
\BDelta^n:=\{ (x_1,\ldots ,x_n)\in I^n|x_1\le\ldots\le x_n \}
\end{equation}
(a function $f:I\to [n]$ such that $f(x)=i$ for
$x_i<x<x_{i+1}$ is identified with 
$(x_1,\ldots ,x_n)\in\Delta^n$; here we assume that
$x_0:=0,x_{n+1}:=n$).

Recall that the nerve $NC$ of a small category $C$ is defined
by $NC(\delta ):=\Funct (\delta ,$ $C)$ for
$\delta\in\Delta_{\byg}$, where $\Funct (\delta ,C)$ is the
set of functors $\delta \to C$. So $|NC|$ is the set of
piecewise constant functors $I\to C$, more precisely 
\begin{equation} \label{nerve}
|NC|:=\rightlimit{F} \Funct (\pi_0(I\setminus F) ),C).
\end{equation}
To pass to the usual definition of $|X|$ use the canonical
representation of a functor $X:\Delta^{\circ}\to\Sets$ as an
inductive limit of representable ones, namely
\begin{equation} \label{limrep}
X=\rightlimit{\Delta /X}\Phi ,
\end{equation}
where $\Delta /X$ is the category of pairs consisting of
$\delta\in\Delta$ and a morphism $\xi :\delta\to X$ (i.e.,
an element $\xi\in X(\delta )$), and $\Phi$ is the functor
from $\Delta /X$ to the category of simplicial sets that
sends $(\delta ,\xi )$ to $\delta$. As geometric realization
commutes with inductive limits and the functor 
$n\mapsto |[n]|$ identifies with the usual ``standard
simplex'' functor $\BDelta :\Delta\to\Sets$
we get the usual formula
\begin{equation} \label{lim}
|X|=\rightlimit{\Delta /X}|\Phi|, \quad 
|\Phi|(\delta ,\xi ):=\mbox{standard simplex }|\delta |,
\end{equation}
which can be rewritten as
\begin{equation} \label{tensor}
|X|=X\times_{\Delta}\BDelta ,
\end{equation}
where 
$X\times_{\Delta}\BDelta$ denotes the coend
of the bifunctor
$X\times\BDelta :\Delta^{\circ}\times\Delta\to\Sets$ 
(which is a kind of ``tensor product of the right
$\Delta$-module $X$ and the left $\Delta$-module $\BDelta$'',
see \S IX.6 of \cite{ML} or Appendix B.7 of \cite{L} for more
details on the $\times_{\Delta}$ bifunctor).

One uses the topology on the standard simplices and either
(\ref{lim}) or (\ref{tensor}) to define the standard topology
on the geometric realization, then one checks that geometric
realization considered as a functor from the
category of simplicial sets to that of compactly generated
spaces commutes with finite projective limits; see
\cite{GZ,GJ,M} for more details.

One can describe the topology on $|X|$ using the
following metric. To define it first choose a measure $\mu$
on $I$ such that the measure of every point is zero and the
measure of every non-empty open set is non-zero. For a
finite $F\subset I$ denote by $\mu_F$ the measure on 
$\pi_0(I\setminus F)$ induced by $\mu$. If
$u,v\in X (\pi_0(I\setminus F))$ define the distance 
$d_{\mu}(u,v)$ to be the minimum of 
$\mu_F (\pi_0(I\setminus F)\setminus A)$ for all subsets
$A\subset\pi_0(I\setminus F)$ such that the images of
$u$ and $v$ in $X(A)$ are equal. If $F'\supset F$ and
$u',v'$ are the images of $u,v$ in $X (\pi_0(I\setminus F))$
then $d_{\mu}(u',v')=d_{\mu}(u,v)$, so we get a well defined
metric $d_{\mu}$ on $|X|$. 

It is easy to see that $d_{\mu}$ is continuous for the
standard topology of $|X|$. It follows that this topology is
Hausdorff. It also follows that if $X$ is finite the
standard topology coincides with the one defined by the
metric $d_{\mu}$. For any $X$ the topological space $|X|$ is
the direct limit of the geometric realizations of the finite
subsets of $X$. 

\medskip

\noindent {\bf Remarks.} (i) Already \cite{Mi} and \cite{GZ}
suggest that in the theory of simplicial sets one should use
the realization (\ref{sim}) of the standard simplex rather
than the ``baricentric'' realization 
$y_0+\ldots y_n=1,y_i\ge 0$. This is also natural in
view of the theory of iterated integrals \cite{H}. 

(ii) Here is a way to think of the  pro-object $\Ifrom$
from (\ref{Ifrom}). Consider the following category $\nabla$.
The objects of $\nabla$ are finite linearly ordered sets $J$
of order $\ge 2$ (so the minimal element 
$\bar 0=\bar 0_J\in J$ does not equal the maximal element 
$\bar 1=\bar 1_J\in J$). The morphisms of $\nabla$ are
nondecreasing maps $f:J\to J'$, $J,J'\in\Ob\nabla$, such
that $f(\bar 0_J)=\bar 0_{J'}$ and 
$f(\bar 1_J)=\bar 1_{J'}$. It is essentially explained in \S
III.1.1 of \cite{GZ} that $\nabla$ is antiequivalent to
$\Delta$. More precisely, the functor 
$\Delta_{\byg}^{\circ}\to \nabla$ that sends
$I\in\Delta_{\byg}$ to the set
$I^*:=\Hom_{\Delta_{\byg}}(I,\{ 0,1\})$ equipped with the
natural (``pointwise'') order is an equivalence. To see this
first notice that it has a left adjoint (namely the functor
$\nabla\to\Delta_{\byg}^{\circ}$ that sends
$J\in\nabla$ to the set $J^*:=\Hom_{\nabla}(J,\{0 ,1\})$
equipped with the natural order) and then show that for all
$I\in\Ob\Delta_{\byg}$ and $J\in\Ob\nabla$ the adjunction
morphisms $I\to I^{**}$ and $J\to J^{**}$ are isomorphisms.

The functor $I\mapsto I^*$
induces an antiequivalence between the category of
pro-objects of $\Delta_{\byg}$ and that of ind-objects of
$\nabla$. The latter category identifies with the category
$\nabla_{\infty}$ of all linearly ordered sets $J$
having a minimal element $\bar 0=\bar 0_J\in
J$ and a maximal element $\bar 1=\bar 1_J\in J$ such that
$\bar 0\ne \bar 1$ (the morphisms of this category are
nondecreasing maps such that $\bar 0\mapsto\bar 0$ and $\bar
1\mapsto\bar 1$). So one can consider $I=[0,1]$ as an
ind-object of $\nabla$. It is easy to see that {\it the
pro-object $\Ifrom$ from (\ref{Ifrom}) corresponds to
$I=[0,1]$ viewed as an ind-object of $\nabla$.} Therefore
(\ref{Ifrom}) implies Theorem III.1.3 of \cite{GZ}, which
says that the group of automorphisms of the geometric
realization functor from the category of simplicial sets to
that of sets equals the group of automorphisms of the
ordered set $I$ (which is the same as the group of
orientation preserving homeomorphisms of $I$). 
 
(iii) Formula (\ref{Ifrom}) is natural in
view of the general theory of fiber functors on a topos
developed in \S 6.8 of \cite{GV} (especially cf. \S 6.8.6 -
\S 6.8.7). We need only the case of a topos of the form
$\hat C:=$the category of presheaves of sets on a small
category $C$  (i.e., $\hat C$ is the category of functors
$C^{\circ}\to\Sets$). According to \cite{GV}, a  functor
$\hat C\to\Sets$ is said to be a {\it fiber functor\,} if it
commutes with finite projective limits and arbitrary
inductive limits. Let $\Fib (\hat C)$ denote the full
subcategory of the category of functors $\hat C\to\Sets$
that consists of fiber functors. It is closed under
inductive limits, so the functor
$F:C^{\circ}\to\Fib (\hat C)$ defined by $(Fc)(X):=X(c)$,
$c\in C$, $X\in\hat C$, canonically extends to a functor
$(\Pro C)^{\circ }=\Ind C^{\circ }\to\Fib (\hat C)$, where
$\Pro C$ (resp. $\Ind C$) denotes the category of
pro-objects (resp. ind-objects). According to \S 6.8 of
\cite{GV}, {\it the functor 
$(\Pro C)^{\circ }\to\Fib (\hat C)$ is an equivalence}.
As explained in \cite{GV}, the quasi-inverse functor is
constructed as follows. Let $\Phi :\hat C\to$Sets be a fiber
functor. Let $\cI$ be the category of pairs $(c,\xi)$, 
$c\in C^{\circ}$, $\xi\in\Phi (c)$. According to \cite{GV},
$\cI$ is a filtering category, and the functor
$\cI\to C^{\circ}$ defines the desired pro-object of $C$. To
show that $\cI$ is filtering the authors of
\cite{GV} first notice that the category $\hat\cI$ of pairs
$(c,\xi)$,
$c\in\hat C^{\circ}$,
$\xi\in\Phi (c)$, is filtering (because $\Phi$ commutes with
finite projective limits) and then deduce from this that 
$\cI$ is also filtering using the fact that $\Phi$
commutes with direct limits (this fact and the analog of
formula (\ref{limrep}) with $\Delta$ replaced by $C$ imply
that every object of $\hat\cI$ can be mapped to an object of
$\cI$).

\section{$\BZ_+$-categories and the category $\Lambda$}  
\label{2}

In the theory of cyclic sets the role of $I=[0,1]$ and
$[n]\in\Delta$ is played by certain $\BZ_+$-categories.

We define a {\it $\BZ_+$-category\,} to be a category $C$
equipped with a morphism of monoids $\BZ_+\to\End\,\id_C$
(the operation in $\BZ_+$ is addition). In other words, each
object $c\in C$ should be equipped with an endomorphism
$\bone_c$ so that $f\bone_{c_1}=\bone_{c_2}f$ for every
morphism $f:c_1\to c_2$ in $C$. If $C,C'$ are
$\BZ_+$-categories then a {\it $\BZ_+$-functor\,} $C\to C'$
is a functor $\Phi:C\to C'$ such that
$\Phi (\bone_{c})=\bone_{\Phi (c)}$ for all $c\in C$. A {\it
$\BZ_+$-isomorphism\,} is a $\BZ_+$-functor which is an
isomorphism. A full subcategory of a $\BZ_+$-category is a
$\BZ_+$-category. 

\medskip

\noindent {\bf Examples.} In this work we will work only
with the following examples of $\BZ_+$-categories. All of
them belong to the class described in Proposition
\ref{equiv}.

 1. {\it Unless stated otherwise,
we will always equip $S^1:=\BR /\BZ$ with the following
structure of $\BZ_+$-category:\,} morphisms from $c\in S^1$
to $c'\in S^1$ are homotopy classes of ``oriented'' paths
$I\to S^1$ (i.e., of those paths which can be represented as
a composition $I\buildrel{f}\over{\to} \BR\to\BR /\BZ=S^1$
with $f$ nondecreasing), and the composition of morphisms is
usual (so our category is a subcategory of the fundamental
groupoid of $S^1$, and $\bone_{x}$ is the homotopy class of a
degree $1$ path from $x\in S^1$ to $x$). 

The category $S^1$ is generated by morphisms $\xy$, $x\ne y$,
where $\xy$ is the homotopy class of the shortest oriented
path from $x\in S^1$ to $y\in S^1$. If $x\ne y$ then
$\yx\xy=\bone_{x}$, and if $x,y,z\in S^1$ are distinct points
in the correct cyclic order then $\yz\xy=\xz$. This is a
complete set of relations between the generators $\xy$ and
$\bone_{x}$. 

2. For every $n\in\BZ_+$ embed $[n]:=\{ 0,\ldots ,n\}$ into
$S^1$ by $k\mapsto k/(n+1) \mod\BZ$. The set $[n]$ equipped
with the $\BZ_+$-category structure induced from $S^1$ will
be denoted by $[n]_{\cyc}$.

3. For every non-empty finite $F\subset S^1$ there is a
unique $\BZ_+$-category $C$ equipped with a $\BZ_+$-functor
$\pi_F :S^1\setminus F\to C$ whose restriction to any system
of representatives $R$ of connected components of 
$S^1\setminus F$ is an isomorphism ($R$ is considered as a
full subcategory of $S^1$). We will denote
$C$ simply by $\pi_0(S^1\setminus F)$. Of course, 
$\pi_0(S^1\setminus F)$ is
(non-canonically) $\BZ_+$-isomorphic to $[n]_{\cyc}$,
$n=\Card F -1$. If $F\subset F'$ there is a unique
$\BZ_+$-functor 
$\Phi: \pi_0(S^1\setminus F')\to\pi_0(S^1\setminus F)$ such
that $\Phi\pi_{F'}=\pi_F$.

4. Given a small category $\cA$ we will define a 
$\BZ_+$-category $\cA_{\cyc}$ with the same set of objects
(in fact, we will use this construction only for categories
that come from linearly ordered sets, and it is not clear if
it is reasonable in general). By definition, $\cA_{\cyc}$ is
generated by $\cA$ and new generators $f^*$ corresponding to
$f\in\Mor\cA$. More precisely, if $f$ is a morphism 
$c_1\to c_2$ then $f^*$ is a morphism $c_2\to c_1$
(intuitively, $f^*=\bone f^{-1}=f^{-1}\bone$). Defining
relations in $\cA_{\cyc}$: $(gf)^*g=f^*$, $f(gf)^*=g^*$.
$\cA_{\cyc}$ becomes a $\BZ_+$-category if one defines 
$\bone_c:c\to c$ by $\bone_c:=({\id_c})^*$. To see this
notice that for every $\cA$-morphism $f:c_1\to c_2$ one has
$f\bone_{c_1}=ff^*f=\bone_{c_2}f$ and
$f^*\bone_{c_2}=f^*ff^*=\bone_{c_1}f^*$ (because
$f^*f=(f\id_{c_1})^*f=(\id_{c_1})^*=\bone_{c_1}$ and
$ff^*=f(\id_{c_2}f)^*=(\id_{c_2})^*=\bone_{c_2}$).
Clearly $\cA\mapsto\cA_{\cyc}$ is a functor from the
category of small categories to that of small
$\BZ_+$-categories.

If $\cA$ is the ordered set $[n]$, $n\in\BZ_+$, considered
as a category then $\cA_{\cyc}$ identifies with the
$\BZ_+$-category
$[n]_{\cyc}$ from Example 2. If $F\subset I:=[0,1]$ is a
finite set containg $\{ 0,1\}$ and $\cA$ is the ordered set
$\pi_0(I\setminus F)$ considered as a category then
$\cA_{\cyc}$ identifies with the $\BZ_+$-category 
$\pi_0(S^1\setminus \bar F)$ from Example 3, where 
$\bar F\subset S^1=\BR/\BZ$ is the image of $F$. If $\cA$ is
the ordered set $[0,1)$ then $\cA_{\cyc}$ identifies with
$S^1$.

\begin{pr}  \label{equiv}
The following properties of a small $\BZ_+$-category $C$ are
equivalent:

(i) $C$ is $\BZ_+$-isomorphic to $\cA_{\cyc}$ for some
linearly ordered set $\cA$;

(ii) every non-empty finite full subcategory of $C$ is
$\BZ_+$-isomorphic to $[n]_{\cyc}$ for some $n\in\BZ_+$.

(iii) every full subcategory of $C$ with $n$ objects,
$n\in\{ 1,2\}$, is $\BZ_+$-isomor\-phic to $[n-1]_{\cyc}$.
\end{pr}

\begin{proof}
It suffices to show that (iii)$\Rightarrow$(i). For every
$x,y\in C$ there exists $f:x\to y$ such that every $g:x\to y$
is uniquely representable as $g=f\bone_x^n$,
$n\in\BZ_+$. This $f$ (which is clearly unique)
will be denoted by $f_{xy}$. Now fix $c\in C$. For every
$x,y\in C$ one has $f_{xy}f_{cx}=f_{cy}\bone_c^m$ for some
$m=m(x,y)\in\BZ_+$. Clearly $m(x,y)+m(y,z)=m(x,z)$ and
$m(x,y)+m(y,x)=1$. Let $\cA\subset C$ be the subcategory with
$\Ob\cA=\Ob C$ whose morphisms are the $f_{xy}$'s
corresponding to those $x,y$ for which $m(x,y)=0$. Then
$\cA$ is a linearly ordered set and $\cA_{\cyc}=C$. 
\end{proof}

Let $\Lambda_{\hug}$ be the category whose objects are small
$\BZ_+$-categories satisfying the conditions of Proposition
\ref{equiv} and whose morphisms are $\BZ_+$-functors. Let
$\Lambda$ (resp.$\Lambda_{\byg}$) be
the full subcategory of $\Lambda_{\hug}$ formed by the
$\BZ_+$-categories $[n]_{\cyc}$, $n\in\BZ_+$ (resp. by
$\BZ_+$-categories which are $\BZ_+$-isomorphic to
$[n]_{\cyc}$ for some $n\in\BZ_+$). 

The following remarks are used only at the end of \S\ref{3}.

\medskip

\noindent {\bf Remarks.} (i) A.~Connes \cite{C} showed that
$\Lambda^{\circ}$ is equivalent to $\Lambda$ (see also
Proposition 6.1.11 of \cite{L}). A more understandable
proof of this equivalence was given by Elmendorf \cite{E}.
Here is a modification of it based on an idea of D.~Arinkin.
If $C_1,C_2$ are $\BZ_+$-categories and $C_1$ is small then
$\ZFunct (C_1,C_2)$ (i.e., the full subcategory of
$\BZ_+$-functors in the category of functors
$C_1\to C_2$) has an obvious structure of $\BZ_+$-category.
So we get a functor
$(\Lambda_{\hug})^{\circ}\to\Lambda_{\hug}$ defined by
$C\mapsto C^*:=\ZFunct (C,[0]_{\cyc})$, where $[0]_{\cyc}$
is the $\BZ_+$-category from Example 2 (it has a single
object $0$, and $\End 0=\BZ_+$). As
$([0]_{\cyc})^{\circ}=[0]_{\cyc}$ one  has 
$\ZFunct (C_1,C_2^*)=(\ZFunct (C_2,C_1))^{\circ}$, so one
gets a canonical $\BZ_+$-functor $F_C:C\to C^{**}$. Finally,
if $C=[n]_{\cyc}$ then $C^*$ is $\BZ_+$-isomorphic to
$[n]_{\cyc}$ and $F_C$ is an isomorphism. So the functor
$(\Lambda_{\hug})^{\circ}\to\Lambda_{\hug}$ induces an
equivalence $(\Lambda_{\byg})^{\circ}\to\Lambda_{\byg}$.

(ii) It is easy to see that in the situation of Example 3 the
$\BZ_+$-category $F^*$ defined in Remark (i) canonically
identifies with $\pi_0(S^1\setminus F)$. Here is an abstract
way to identify them. Let $\Ztors$ denote the
$\BZ_+$-category of $\BZ_+$-torsors, i.e., $\BZ_+$-sets
isomorphic to $\BZ_+$. For $C\in\Lambda_{\hug}$ one has the
$\BZ_+$-bifunctor $C^{\circ}\times C\to\Ztors$ defined by
$(c,c')\mapsto\Mor (c,c')$. Composing it with the unique
$\BZ_+$-functor $\Ztors\to [0]_{\cyc}$ one gets a
$\BZ_+$-functor $\Phi=\Phi_C:C\to C^*$. One also has the
$\BZ_+$-functor 
$\Phi'=(\Phi_{C^{\circ}})^{\circ}:C\to
((C^{\circ})^*)^{\circ}=C^*$ (which is usually different from
$\Phi_C$). In the situation of Example 3 the compositions
$S^1\setminus F\mono S^1
\buildrel{\Phi}\over{\longrightarrow}(S^1)^*\to F^*$ and
$S^1\setminus F\mono S^1
\buildrel{\Phi'}\over{\longrightarrow}(S^1)^*\to F^*$ are
equal. This $\BZ_+$-functor can be taken as $\pi_F$ (see
Example 3).

(iii) The category of ind-objects of
$\Lambda_{\byg}$ (or of $\Lambda$) identifies with
$\Lambda_{\hug}$. The category of pro-objects of
$\Lambda_{\byg}$ can be identified with 
$(\Lambda_{\hug})^{\circ}$ using the equivalence
$(\Lambda_{\byg})^{\circ}\iso\Lambda_{\byg}$ from Remark (i).
It gives a class of equivalences $\Lambda^{\circ}\iso\Lambda$
each two of which are related by a unique isomorphism.

\section{Cyclic sets}     \label{3}

In \S\ref{2} we defined the categories
$\Lambda_{\hug}\supset\Lambda_{\byg}\supset\Lambda$.
According to A.~Connes, a {\it cyclic set\,} is a functor
$X:\Lambda^{\circ}\to\Sets$. One can extend such $X$ to a
functor $\Lambda_{\byg}^\circ\to\Sets$, and the
extension is unique up to unique isomorphism. The functor
$C\mapsto C_{\cyc}$ from Example 4 of \S\ref{2} can be
considered as a functor $\Delta_{\byg}\to\Lambda_{\byg}$.
Using this functor one can consider any cyclic set as a
simplicial set.

The {\it geometric realization} $|X|$ of a cyclic set $X$
is defined in \cite{GZ,GJ,M} to be the geometric
realization of $X$ considered as a simplicial set. By
(\ref{defcyc}), this means that $|X|$ is the direct limit of
$X((\pi_0(I\setminus F))_{\cyc})$, where $F$ runs through the
set of all finite subsets of
$I:=[0,1]$. We will get the same answer if $F$ runs
only through the set of finite subsets of $I$ containing $0$
and $1$. For such $F$ one has 
$(\pi_0(I\setminus F))_{\cyc}=\pi_0(S^1\setminus\bar F)$,
where $\bar F\subset S^1=\BR/\BZ$ is the image of $F$ (see
Example 4 of \S\ref{2}). So 
\begin{equation} \label{defcyc}
|X|:=\rightlimit{F} X (\pi_0(S^1\setminus F) ),
\end{equation}
where $F$ runs through the set of all non-empty finite
subsets of $S^1$ and $\pi_0(S^1\setminus F)$ is equipped with
the $\BZ_+$-category structure from Example 3 of \S\ref{2}. 
The reader may prefer to rewrite (\ref{defcyc}) as 
\begin{equation} \label{Sfrom}
|X|:=X(\Sfrom ), 
\end{equation}
where $\Sfrom$ is the pro-object of
$\Lambda_{\byg}$ which is the projective limit of the objects
$\pi_0(S^1\setminus F)\in\Lambda_{\byg}$ over all finite
subsets $F\subset S^1$. 

Let us recall that the cyclic nerve $\N C$ of a small
category $C$ is defined by 
$\N C(\lambda )=\Funct (\lambda ,C)$ for
$\lambda\in\Lambda_{\byg}$, where $\Funct (\lambda ,C)$ is
the set of functors $\lambda\to C$. So $|\N C|$ is the set of
piecewise constant functors $S^1\to C$, more precisely 
\begin{equation} \label{nervecyc}
|\N C|:=\rightlimit{F} \Funct (\pi_0(S^1\setminus F) ),C).
\end{equation}
   
If $X$ is a cyclic set then $|X|$ is the geometric
realization of $X$ considered as a simplicial set, so it
is equipped with a topology.    On the other hand, formula (\ref{defcyc}) or
(\ref{Sfrom}) makes it clear that the group $\Aut S^1$ of
orientation preserving homeomorphisms of $S^1$ acts on
$|X|$. In particular, the group of rotations $SO(2)$ acts. 

We claim that the action of $\Aut S^1$ on
$|X|$ is continuous if $\Aut S^1$ is equipped with the
compact-open topology. To prove this, note
that similarly to (\ref{tensor}) and
(\ref{lim}) one has canonical homeomorphisms
\begin{equation} \label{tensorcyc}
|X|=X\times_{\Lambda}\Psi,
\end{equation}
\begin{equation} \label{limcyc}
|X|=\rightlimit{\Lambda/X}\Psi_X,
\end{equation}
where $\Psi :\Lambda\to\mbox{\{Spaces\}}$, $\Psi (\lambda ):=$the
geometric realization of the functor represented by
$\lambda$, and $\Psi_X$ is the restriction of $\Psi$ to
$\Lambda/X$. Using the homeomorphism \eqref{tensorcyc} or  \eqref{limcyc} one reduces the question to the particular case where $X$ is the cyclic set $[n]_{\cyc}$ 
(i.e., the functor $\Lambda^{\circ}\to\Sets$ represented by $[n]_{\cyc}$). 
In this case one can use the following explicit description of the geometric realization $|[n]_{\cyc}|$.

\medskip

\noindent {\bf The geometric realization of the cyclic set $[n]_{\cyc}$.} 
We claim that there are canonical bijections
\begin{equation}   \label{e:bijections}
|[n]_{\cyc}|\iso\Funct_{\BZ_+}([n]_{\cyc},S^1)\iso Y/\BZ,
\end{equation}
where $\Funct_{\BZ_+}([n]_{\cyc},S^1)$ stands for the set of $\BZ_+$-functors $[n]_{\cyc}\to S^1$ and
$Y$ is the set of non-decreasing maps $f:\BZ\to\BR$ such that $f(i+n+1)=f(i)$ for all $i\in\BZ$ (the $\BZ$-action on $Y$ is induced by the action of $\BZ$ on $\BR$ by translations). The second bijection in \eqref{e:bijections} is clear. The first one is constructed as follows. By \eqref{defcyc},
\[
|[n]_{\cyc}|:=\rightlimit{F}  \Funct_{\BZ_+}(\pi_0(S^1\setminus F),[n]_{\cyc}).
\]
By Remark~(ii) from \S\ref{2}, $\pi_0(S^1\setminus F)=F^*$, so after fixing a $\BZ_+$-isomorphism
$([n]_{\cyc})^*\iso [n]_{\cyc}\,$, the set  $\Funct_{\BZ_+}(\pi_0(S^1\setminus F),[n]_{\cyc})$ identifies with
\[
\Funct_{\BZ_+}([n]_{\cyc}\, , F)
\]
and the set $|[n]_{\cyc}|$ identifies with
\[
\rightlimit{F} \Funct_{\BZ_+}([n]_{\cyc}\, , F)= \Funct_{\BZ_+}([n]_{\cyc}\, , S^1).
\]
It is easy to check that the composition of the bijections \eqref{e:bijections} is a homeomorphism
(if $Y$ is equipped with the natural topology). Note that $$\Aut S^1=G/\BZ,$$ where $G$ is the group of orientation-preserving homeomorphisms $f:\BR\iso\BR$ such that $f(x+1)=f(x)+1$. So $\Aut S^1$ acts on 
$Y/\BZ$. This action is continuous, and the homeomorphism $|[n]_{\cyc}|\iso Y/\BZ$ is ($\Aut S^1$)-equivariant.
So the action of $\Aut S^1$ on $|[n]_{\cyc}|$ is continuous.


\medskip
\noindent {\bf Remark.} 
The space $|[n]_{\cyc}|=Y/\BZ$ is homeomorphic to the product of $S^1$ and the $n$-dimensional simplex.

\medskip

\noindent {\bf Warning.} 
Associating to a $\BZ_+$-functor $\Phi : [n]_{\cyc}\to S^1$ the sequence $\Phi (0),\ldots,\Phi (n)$
one gets a map
\[
|[n]_{\cyc}|=\Funct_{\BZ_+}([n]_{\cyc},S^1)\to (S^1)^{n+1}.
\]
Unless $n=0$, it is \emph{not injective} (contrary to what I wrote in the previous versions of this article\footnote{The mistake was noticed by M.~Kapranov.}, including the journal version). More precisely, the preimage of a point of the diagonally embedded 
$S^1\subset (S^1)^{n+1}$ has cardinality $n+1$. (All other preimages have cardinalities $\le 1$.)

\medskip

\noindent {\bf Remarks.} (i) By Remark (iii) at the end of
\S\ref{2}, the pro-object $\Sfrom$ from (\ref{Sfrom}) can be
considered as an object of $\Lambda_{\hug}$. This object
canonically identifies with $S^1$ (e.g., one can use Remark
(iii) from \S\ref{2}).

(ii) The group of automorphisms of the geometric
realization functor from the category of cyclic sets to
that of sets equals $\Aut S^1$. This follows, e.g., from the
previous remark.

\end{document}